\documentclass[final,leqno,onefignum,onetabnum]{siamltex1213}
\usepackage[T1]{fontenc}
\usepackage{amsfonts}
\usepackage{amsmath}
\usepackage{amssymb}
\usepackage{bbm}
\usepackage{bm}
\usepackage{mathrsfs}
\usepackage{verbatim}
\usepackage{setspace}
\usepackage{color}
\usepackage{pdfsync}
\usepackage{enumitem}
\usepackage{scrtime}
\newtheorem{remark}[theorem]{Remark}

\newcommand*{\qedh}{\hfill\ensuremath{\square}}%

\def \sT {s\in [t,T]}
\def \rw {\rightarrow}

\def \tx {(t,x)\in [0,T]\times \bR}
\newcommand{\cF}{\mathcal{F}}
\newcommand{\cP}{\mathcal{P}}
\newcommand{\cH}{\mathcal{H}}
\newcommand{\cM}{\mathcal{M}}
\newcommand{\bP}{\textbf{P}}
\newcommand{\E}{\textbf{E}}
\newcommand{\bR}{\textbf{R}}
\newcommand{\tT}{t\leq T}
\newcommand{\no}{\noindent}
\newcommand{\su}{u^*}
\newcommand{\sv}{v^*}
\newcommand{\bu}{\bar{u}}
\newcommand{\bv}{\bar{v}}
\newcommand{\Y}{Y^{u,v}}
\newcommand{\Z}{Z^{u,v}}

\newcommand{\be}{\begin{equation}}
\newcommand{\ee}{\end{equation}}
\newcommand{\ba}{\begin{array}}
\newcommand{\ea}{\end{array}}

\def \ms {\medskip}
\def \G {\Gamma}
\def \lb {\label}
\def \lg {\lambda}
\def \sp{[0,T]\times \bR}

\title{On the Bang-Bang Type Nash Equilibrium Point for Markovian Nonzero-sum Stochastic Differential Game} 

\author{
  Said Hamad\`ene
 \thanks{
 Universit\'e du Maine, LMM, Avenue Olivier Messiaen, 72085 Le Mans, Cedex 9, France. (\email{hamadene@univ-lemans.fr}).
 }
  \and Rui Mu \thanks{ Corresponding author. Center for Financial Engineering, Soochow University, Suzhou 215006, P. R. China;
  Universit\'e du Maine, LMM, Avenue Olivier Messiaen, 72085 Le Mans, Cedex 9, France. (\email{rui.mu.sdu@gmail.com}). This author's research was supported in part by the Natural Science Foundation for Young Scientists of Jiangsu Province, P. R. China (No. BK20140299).}
 }

\begin{document}
\maketitle
\slugger{mms}{xxxx}{xx}{x}{x--x}

\begin{abstract}
In this paper, we study a nonzero-sum stochastic differential game of bang-bang type in the Markovian framework. We show the existence of a Nash equilibrium point for this game. The main tool is the notion of backward stochastic differential equations which, in our case, are multidimensional with discontinuous generators with respect to $z$ component.
\end{abstract}



\pagestyle{myheadings}
\thispagestyle{plain}
\markboth{TEX PRODUCTION}{USING SIAM'S \LaTeX\ MACROS}

\section{Introduction}\label{se:intro}
This paper is related to nonzero-sum stochastic differential games (NZSDG, for short) of bang-bang type in the Markovian framework which we describe below. 

We consider the case of two players $\pi_1$ and $\pi_2$. If there are more than two players, the adaptation is straightforward. The two players $\pi_1$ and $\pi_2$ act on a system through two admissible controls $u:=(u_s)_{s\leq T}$ and $v:=(v_s)_{s\leq T}$ which are adapted stochastic processes, with values in compact metric spaces $U$ and $V$, respectively. The dynamics of the controlled system is given by a stochastic process $(X^{u,v}_s)_{s\leq T}$, solution of the following standard stochastic differential equation (SDE, for short):
\begin{equation}\label{eq: intro x}
dX_s^{u,v}=\Gamma(s,X_s^{u,v},u_s,v_s)ds+\sigma(s,X_s^{u,v})dB_s\text{ for }s\leq T, \text{ and } X_0=x_0
\end{equation}
where $B:=(B_s)_{s\leq T}$ is a Brownian motion. Next with each player $\pi_i$, $i=1,2$, is associated a payoff $J_i(u,v)$, $i=1,2$, given by:
\begin{equation}\label{eq: payoff intro}
J_i(u,v)=\E[g_i(X_T^{u,v})].
\end{equation}
The objective is to find a pair $(u^*,v^*)$ which satisfy 
\begin{equation}\label{eq: obj intro}
J_1(u^*, v^*)\geq J_1(u, v^*) \text{ and } J_2(u^*, v^*)\geq J_2(u^*, v)\mbox{ for any pair }(u,v).
\end{equation}

As we can see, the payoff function of the player $\pi_1$ (resp. $\pi_2$) depends not only on its own control $u$ (resp. $v$) but also on the control used by the other player $\pi_2$ (resp. $\pi_1$). Therefore there is a game  between the two players which is of cooperative relationship. This kind of game is known as the nonzero-sum stochastic differential game. The pair $(u^*,v^*)$ of (\ref{eq: obj intro}) is called a Nash equilibrium point (NEP, for short) of the game. Its meaning is that when it is implemented by both players, then if $\pi_1$ (resp. $\pi_2$) decides unilaterally to change a control while $\pi_2$ (resp. $\pi_1$) keeps $v^*$ (resp. $u^*$) then her payoff $J_1(u, v^*)$ (resp. $J_2(u^*, v)$) is lesser than 
$J_1(u^*, v^*)$ (resp. $J_2(u^*, v^*)$), i.e., her action of unilateral deviation induces a penalty.

In the case when $J_1(u,v)+J_2(u,v)=0$, the game turns into the well-known zero-sum differential game which is widely studied in the literature (see e.g. \cite{zs3, zs4, zs5, zs1, zs2}, etc. and the references therein). On the other hand, if $X^{u,v}$ does not depend on $v$ then the problem turns merely into  a control problem. In this specific case, we know that an optimal control exists and is of bang-bang type since it takes values only on the boundaries of $U$ according to the derivative of the value function of the control problem. This is the consequence of the fact that the instantaneous reward in (\ref{eq: payoff intro}) is null. So one would expect the same features of the NEP of this game if it exists.
\ms

The nonzero-sum differential game is also considered by several authors in the literature, see eg. \cite{bensoussan1, bensoussan2,buckdahn2004,  carmona2013control, rainer, friedman2, H1997, hamadene1999, hamadene1998, qianlin, mannucci, carmona2013probabilistic}, to name a few. There are typically two approaches. One method is related to partial differential equation (PDE) theory. Some of the results show that the payoff function of the game is the unique viscosity solution of a related Hamilton-Jacobi-Bellman-Isaacs equation, e.g., \cite{qianlin}. Other works make use of the Sobolev theory of PDEs (see \cite{bensoussan1, bensoussan2, friedman2, mannucci}, etc.) to deal with NZSDGs. Comparatively, another popular way to deal with stochastic differential game is the backward stochastic differential equation (BSDE) approach \cite{H1997, hamadene1999, hamadene1998, qianlin}, which characterizes the payoffs of the game through  solutions of associated BSDEs. However those BSDEs are of multidimensional type and usually their generators are non-Lipschitz. Therefore proving that they have solutions is not an easy task. 
\medskip

In the present article, we study the above NZSDG via the BSDE arguments in the Markovian framework. For clarity reasons we consider a special game model in assuming that: 

(i) the process $X$ of \eqref{eq: intro x} is $\bR$-valued,  and $U=[0,1]$, $V=[-1,1]$ ;

(ii) the drift $\Gamma$ of \eqref{eq: intro x} has the following structure:
$$
\Gamma(t,x,u,v)=f(t,x)+u+v.
$$
The conditions on the functions $f$ and $g_i$, $i=1,2$, are rather weak since they are related to measurability and growth conditions. 

At the end of the paper we give hints which allow to generalize this setting in several directions, especially if the dynamics contains a diffusion term and the multidimensional case for the process $X^{u,v}$. 
\ms

In this problem the main difficulty is lodged at the level of the main BSDE (\ref{eq:main bsde}) associated with this NZSDG. Once the existence of its solution is stated, it provides the NEP of the game. As pointed out previously, this BSDE is of multidimensional type (here of dimension two since there are two players) and whose generator is discontinuous in $z$. The main challenge we overcome is to show that this BSDE has a solution and then we constructed a NEP for the game defined by (\ref{eq: intro x}), (\ref{eq: payoff intro}) and (i)-(ii) above. Like in the control framework, this NEP is of bang-bang type
since the payoffs have no instantaneous payoffs included. This is the main novelty of this article. The closest work to ours is the one by P.Manucci \cite{mannucci}, but this latter concerns only diffusions in bounded domains and the requirements on the regularity of the data are stronger than ours due to the method she employed based on PDEs in Sobolev spaces.  
\ms

This paper is organized as follows:
\medskip

In Subsection \ref{subsec:statement 1-d}, we introduce the game problem and some preliminaries. The formulation we adopt is of weak type. Besides, for intuitive understanding, we work with a particular setting of controls and state process $X^{u,v}$. The extension to the multidimensional situation obviously holds following the same ideas. The explicit form of discontinuous controls, namely, bang-bang controls are presented in Subsection \ref{subsec: bang-bang}. In Subsection \ref{subsec: main result}, we give the main result (Theorem \ref{th: nash}) of this work and some other related important results. We first provide a link between the game problem and Backward SDEs (Proposition \ref{prop: yuv}). The payoff of the game can be characterized by the initial value of the solution for an associated BSDE. Then, by Proposition \ref{prop: equilibrium}, we prove that the existence of a NEP for the game is equivalent to the existence of a solution of a BSDE which is of multidimensional and with discontinuous generator with respect to $z$. Finally, under some reasonable assumption, we provide the solution of this special BSDE (Theorem \ref{th: main bsde}). All the proofs are stated in Subsection \ref{subsec:proofs}. The proofs of Propositions \ref{prop: yuv} and \ref{prop: equilibrium} are standard. For Theorem \ref{th: main bsde}, the method is mainly based on an approximating scheme. In Section \ref{sec:gener}, we investigate some possible generalizations. The idea is the same with a bit modification which is indicated. 

\section{Bang-bang type NZSDG and multidimensional BSDEs with discontinuous generators}\label{sec: 1-d}

In this section, we first deal with the bang-bang type nonzero-sum stochastic differential game problem in 1-dimensional framework. A more general setting will be given in the next section. 
\subsection{Statement of the problem}\label{subsec:statement 1-d}

Let $T>0$ be fixed and let $(\Omega, \cF, \bP)$ be a probability space, on which is defined a 1-dimensional standard Brownian motion $B:=(B_t)_{t\leq T}$. For $t\leq T$, let us set $F_t:=\sigma(B_u,u\leq t)$ and denote by $(\cF_t)_{t\leq T}$ the completion of $(F_t)_{t\leq T}$ with the $\bP$-null sets of $\cF$. Next let $\cP$ be the $\sigma$-algebra on $[0,T]\times \Omega$ of ${\cF}_t$-progressively measurable sets. For a real constant $p\geq 1$, we introduce the following spaces:
\no \begin{itemize}
\item $L^p=\{\xi:\ \cF_T$-measurable and $\bR$-valued random variable s.t.$\E[|\xi|^p]<\infty\}$;
\item $\mathcal{S}_T^p(\bR)=\{Y=(Y_s)_{s\in[0,T]}:\ \mathcal{P}$-measurable, continuous and $\bR$-valued stochastic process s.t. $\E[\sup_{s\in [0,T]}|Y_s|^p]<\infty\}$;
\item $\mathcal{H}_T^p(\bR)=\{Z=(Z_s)_{s\in[0,T]}:\ \mathcal{P}$-measurable and $\bR$-valued stochastic process s.t. $\E[(\int_0^T|Z_s|^2ds)^{p/2}]<\infty\}$.
\end{itemize}

We consider, in this article, the 2-player case which we describe accurately below. The general multiple players case is a straightforward adaptation.

Let $(t,x)\in [0,T]\times \bR$ and $X^{t,x}$ be the stochastic process defined as follows:
\begin{equation}\label{eq: SDE}
\forall s\leq T,\ X_s^{t,x}=x+(B_{s\vee t}-B_t).
\end{equation}

\begin{remark}\label{re:sigma} Note that we consider a trivial situation for SDE \eqref{eq: SDE} with an identity diffusion process, just for easy understanding. The trick of the technique in this article still valid for general diffusion process with appropriate properties. We will introduce this point in Section \ref{sec:gener}.
\end{remark}

Each player $\pi_i$, $i=1, 2$, has her own control. Let us denote next by $U=[0,1]$, $V=[-1,1]$ those two compact subsets of $\bR$ and $\cM_1$ (resp. $\cM_2$) the set of $\cP$-measurable process $u=(u_t)_{\tT}$ (resp. $v=(v_t)_{\tT}$) on $[0,T]\times \Omega$ with value in $U$ (resp. $V$). Hereafter, we call $\cM:=\cM_1\times \cM_2$ (resp. $\cM_1$, resp. $\cM_2$) the set of admissible controls for the two players (resp. first player ; resp. second player). 
\medskip

Let $f: (t,x)\in [0,T]\times \bR\rightarrow \bR$ be a Borelian function.
We will say that $f$ is of linear growth if there exists a constant $C\geq 0$ such that for any $\tx$, 
\begin{equation}\lb {lg}|f(t,x)|\leq C(1+|x|).\end{equation}

Next let $\Gamma$ be the function such that for any $(t,x,u,v)\in [0,T]\times \bR\times U\times V$ associates $\Gamma(t,x,u,v)=f(t,x)+u+v \in \bR$. The function $\Gamma$ stands for the drift of the dynamics of the system when controlled by the two players $\pi_i$, $i=1, 2$. When $f$ is of linear growth, the function $\G$ is so since $U$ and $V$ are bounded sets.
\ms

Next for $(u.,v.)\in \cM$, let $\bP^{u,v}_{t,x}$ be the positive measure on $(\Omega,\cF)$ defined as follows:
\begin{equation}\label{eq:density fun}\begin{array}{c}
d\bP^{u,v}_{t,x}=\zeta_T(\Gamma(.,X.^{t,x},u.,v.))d\bP\ \text{with}\ \zeta_s(\Theta):=1+\int_0^s \Theta_r \zeta_r dB_r,\ \sT,\end{array}
\end{equation}
for any $\cP$-measurable $\bR$-valued process $\Theta:=(\Theta_s)_{s\leq T}$. It follows from the uniform linear growth property of $\Gamma$ that $\bP^{u,v}_{t,x}$ is a probability on $(\Omega, \cF)$ (see Appendix A of \cite{karoui and hamadene} or \cite{karatzas}, pp.200). Then, by Girsanov's Theorem (\cite{girsanov}), the process $B^{u,v}=(B_s-\int_0^s\Gamma(r, X_r^{t,x},u_r,v_r)dr)_{s\leq T}$ is a $(\cF_s, \bP^{u,v}_{t,x})$-Brownian motion and $(X_s^{t,x})_{s\leq T}$ satisfies the following SDE
\begin{equation}\label{eq: SDE with Gamma}\begin{array}{c}
dX_s^{t,x}=\Gamma(s, X_s^{t,x},u_s,v_s)ds+dB_s^{u,v},\ \forall s\in[t,T]\ \text{and}\ X_s^{t,x}=x,\ s\in [0,t].\end{array}
\end{equation}
As a matter of fact, the process $X^{t,x}$ is not adapted with respect to the filtration generated by the Brownian motion $B^{u,v}$. Thereby, $X^{t,x}$ is a weak solution for the SDE \eqref{eq: SDE with Gamma}. If the system starts from $x_0\in \bR$ at $t=0$ and is controlled by player $\pi_1$ (resp. $\pi_2$) with $u.$ (resp. $v.$), the law of its dynamics is the same as the one 
of $X^{0,x_0}$ under $\bP^{u,v}_{0,x_0}$.
\medskip

Once more let $x_0\in \bR$ fixed. We will precise the payoffs of the players when they implement the pair of strategies $(u.,v.)$. It is of terminal type and given, for player $\pi_1$ (resp. $\pi_2$), by 
\begin{eqnarray}\label{eq: payoffs}
J_1(u,v):=\E^{u,v}_{0,x_0}[g_1(X_T^{0,x_0})] 
\mbox{ (resp. }J_2(u,v):=\E^{u,v}_{0,x_0}[g_2(X_T^{0,x_0})]),
\end{eqnarray}
where:

\no (i) $g_1$ and $g_2$ are two Borel measurable functions from $\bR$ to $\bR$ which are of polynomial growth, \textit{i.e.}, there exist non-negative constants $C$ and $\gamma\geq 1$ such that for any $x\in \bR$, 
\begin{eqnarray}\label{(A2)}|g_1(x)|+|g_2(x)|\leq C(1+|x|^{\gamma})\end{eqnarray}
(ii) For any fixed $(t,x)\in \sp$, $\E^{u,v}_{t,x}$ is the expectation under the probability $\bP^{u,v}_{t,x}$ ; hereafter $\E^{u,v}_{0,x}(.)$ (resp. $\bP^{u,v}_{0,x}$) will be simply denoted by $\E^{u,v}_x(.)$ 
(resp. $\bP^{u,v}_x$).
\medskip

As we can see from \eqref{eq: SDE with Gamma} and \eqref{eq: payoffs}, the choice of control of each player has influence on the other one's payoff through the state process $X^{0,x_0}$ under $\bP^{u,v}_{x_0}$. What we discussed here is a nonzero-sum stochastic differential game which means that the two players are of cooperate relationship. Both of them want to reach the maximum payoff. Therefore, naturally, we are concerned with the existence of a \textit{Nash equilibrium point}, which is a couple of controls $(\su,\sv)\in \cM$, such that, for all $(u,v)\in \cM$,
\[
J_1(\su,\sv)\geq J_1(u, \sv)\ \text{and}\ J_2(\su,\sv)\geq J_2(\su,v).
\]
This means that when the strategy $(\su,\sv)$ is implemented by the players, one who makes unilaterally the decision to deviate or to change a strategy, while the other one keeps its own choice, is penalized.\qedh

\subsection{Bang-bang type control}\label{subsec: bang-bang}
As pointed out in \eqref{eq: payoffs}, there are no instantaneous payoffs in $J_1$ and $J_2$. Therefore, in comparison with optimal control which is a particular case of our problem (see e.g. \cite{balak, benes, christopeit, maurerosmo}), the equilibrium point of this game, if exists, should be of bang-bang type, \textit{i.e.}, the optimal control $u^*$ (resp.$v^*$) will jump between the two bounds of the value set $U$ (resp.$V$). 
\medskip

To proceed, let $H_1$ and $H_2$ be the \textit{Hamiltonian functions} of this game problem, \textit{i.e.}, the functions (which do not depend on $\omega$) defined from $[0,T]\times \bR\times \bR\times U\times V$ into $\bR$ by:
\begin{eqnarray*}
H_1(t,x,p,u,v)&:=p\Gamma(t,x,u,v)=p(f(t,x)+u+v);\\
H_2(t,x,q,u,v)&:=q\Gamma(t,x,u,v)=q(f(t,x)+u+v).
\end{eqnarray*}

Next let $\epsilon_1$ and $\epsilon_2$ be two arbitrary elements of $U$ and $V$ respectively. Let $\bar{u}$ and $\bar{v}$ be two functions defined on $\bR\times U$ and $\bR\times V$, valued on $U$ and $V$ respectively, as follows: $\forall p,q\in \bR$,
\begin{equation}
\bar{u}(p,\epsilon_1)=
\left\{
\begin{aligned}
1, \ p>0,\\
\epsilon_1,\ p=0,\\
0,\ p<0,
\end{aligned}
\right.
\quad
\text{and}
\quad
\bar{v}(q,\epsilon_2)=
\left\{
\begin{aligned}
1,\ q>0,\\
\epsilon_2,\ q=0,\\
-1,\ q<0.
\end{aligned}
\right.
\end{equation}
Then, we can easily check that $\bar{u}$ and $\bar{v}$ satisfy the \textit{generalized Isaacs' condition} which reads as follows: $\forall \,\,(t,x,p,q,u,v)\in [0,T]\times \bR\times \bR\times U\times V$, 
\begin{equation}\label{eq:Isaacs}
\left\{
\begin{aligned}
H_1^*(t,x,p,q,\epsilon_2)&:=H_1(t,x,p,\bar{u}(p,\epsilon_1),\bar{v}(q,\epsilon_2))\geq H_1(t,x,p,u,\bar{v}(q,\epsilon_2))\mbox{ and }\\
H_2^*(t,x,p,q,\epsilon_1)&:=H_2(t,x,q,\bar{u}(p,\epsilon_1),\bar{v}(q,\epsilon_2))\geq H_2(t,x,q,\bar{u}(p,\epsilon_1),v).
\end{aligned}
\right.
\end{equation}
\begin{remark}
Let us notice that the function $H_1^*$ (resp. $H_2^*$) does not depend on $\epsilon_1$ (resp. $\epsilon_2$) since, $p\bar{u}(p, \epsilon_1)=p\vee 0$ (resp. $q\bar{v}(q,\epsilon_2)=|q|$) does not depend on $\epsilon_1$ (resp. $\epsilon_2$). Besides, they are discontinuous w.r.t. $(p,q)$ since $\bar{v}$ and $\bar{u}$ are so.
\end{remark}
\medskip

We next give the main result of this article without proofs for intuitive understanding. All the proofs are given in Subsection \ref{subsec:proofs}. 
\subsection{Main result}\label{subsec: main result}

As in several papers on the same subject (\cite{karoui and hamadene, H1997, H2014}, etc.), we will adopt the BSDE approach in order to show that this 
particular nonzero-sum stochastic differential game has a Nash equilibrium point.  For sake of clarity, in this subsection we give the main result and the intermediary ones which we need. We collect all their proofs in the next subsection. To begin with, the following result characterizes the payoffs (\ref{eq: payoffs}) through a solution of a multidimensional BSDE.

\begin{proposition}\label{prop: yuv}
Assume that (\ref{lg}) and (\ref{(A2)}) are satisfied. Then for all $(u,v)\in \cM$ and $i=1,2$, there exists a unique pair of $\cP$-measurable processes\\ $(Y^{i;x_0;u,v}, Z^{i;x_0;u,v})$, with values in $\bR\times \bR$, such that: For $i=1,2$, 
\begin{description}
\item[(i)] for all constant $q\geq 1$, 
\begin{equation}\label{eq:esti yuv and zuv}
\E_{x_0}^{u,v}\Big[\sup_{s\in [0,T]}|Y_s^{i;x_0;u,v}|^q+(\int_0^T |Z_s^{i;x_0;u,v}|^2ds)^{\frac{q}{2}}\Big]<\infty.
\end{equation}
\item[(ii)]
\begin{equation}\label{eq:BSDE yiuv}\left\{\begin{array}{l}
-dY_t^{i;x_0;u,v}=H_i(s, X_s^{0,x_0},Z_s^{i;x_0;u,v},u_s,v_s)ds - Z_s^{i;x_0;u,v}dB_s,\,\,s\leq T; \\ Y_T^{i;x_0;u,v}=g_i(X_T^{0,x_0}).\end{array}\right.
\end{equation}
\item[]Moreover $Y_0^{i;x_0;u,v}=J_i(u,v)$.\qedh
\end{description}
\end{proposition}
\medskip

The following result is a verification theorem for the existence of NEP of the game of bang-bang type. 
\begin{proposition}\label{prop: equilibrium} Assume that (\ref{lg}) and (\ref{(A2)}) are satisfied. 
Besides, suppose that there exist two deterministic functions $\eta^1$, $\eta^2$ and stochastic processes $(Y^1, Z^1)$, $(Y^2, Z^2)$ and $\theta$, $\vartheta$ such that:
\begin{description}

\item[(i) (a)]  $\theta$ (resp. $\vartheta$) is a $\cP$-measurable process with values in $U$ (resp. $V$) and $(Y^1, Z^1)$ and $(Y^2, Z^2)$ are two couples of $\cP$-measurable processes $\bR^{1+1}$-valued which satisfy:
 \item[(b)] for $i=1,2$ $\bP$-a.s. $(Z^i_s(\omega))_{s\leq T}$ is $ds$-square integrable and for all $s\leq T$,
  \begin{equation}\label{eq:main bsde}
  \left\{
     \begin{aligned}
     -dY_s^1&=H_1^*(s, X_s^{0,x_0},Z_s^1, Z_s^2,\vartheta_s)ds- Z_s^1dB_s,\ Y_T^1=g_1(X_T^{0,x_0});\\
     -dY_s^2&=H_2^*(s, X_s^{0,x_0},Z_s^1, Z_s^2,\theta_s)ds- Z_s^2dB_s,\ Y_T^2=g_2(X_T^{0,x_0}).
     \end{aligned}
  \right.
  \end{equation}
\item[(ii)] $\eta^1$ and $\eta^2$ are two deterministic measurable functions with polynomial growth from $[0,T]\times \bR$ to $\bR$ such that 
$\bP$-a.s., $\forall s\leq T$, $Y_s^i=\eta^i(s, X_s^{0,x_0})$.
\end{description}
Then, the pair of controls $(\bar{u}(Z_s^1, \theta_s), \bar{v}(Z_s^2, \vartheta_s))_{s\leq T}$ is a bang-bang type Nash equilibrium point of the nonzero-sum stochastic differential game.\qedh
\end{proposition} 
\medskip

Finally since the diffusion coefficient in equation (\ref{eq: SDE}) is equal to the identity and in using a result by El-Karoui et al. \cite{el karoui} which allows the representation of solutions of BSDEs through deterministic functions, in the markovian case of randomness, we prove the existence of processes and deterministic functions which satisfy the requirements of Proposition \ref{prop: equilibrium}. The main difficulty relies on the discontinuity of the generator $H_1^*$ (resp. $H_2^*$) w.r.t. $(p,q)$ which comes from the discontinuity of $\bar{v}$ (resp. $\bar{u}$) on $q=0$ (resp. $p=0$). However we can overcome this difficulty and we have: 
\begin{theorem}\label{th: main bsde} Assume that 
$f$ and $g_i$, $i=1,2$, 
satisfy to (\ref{lg}) and (\ref{(A2)}) respectively. Then there exist $\eta^1$, $\eta^2$, $(Y^1\!,Z^1\!)$, $(Y^2,Z^2)$ and $\theta$, $\vartheta$ which satisfy (i),(a)-(b) and (ii) of Proposition \ref{prop: equilibrium}.
\end{theorem}

As a consequence of Theorem \ref{th: main bsde} and Proposition \ref{prop: equilibrium}, we obtain the main result of this article.
\bigskip

\begin{theorem}\label{th: nash} The pair $(\bar{u}(Z_s^1, \theta_s), \bar{v}(Z_s^2, \vartheta_s))_{s\leq T}$ of $\cM$ is a bang-bang type Nash equilibrium point for the nonzero-sum stochastic differential game defined by (\ref{eq: SDE}), (\ref{eq:density fun}) and (\ref{eq: payoffs}).  

\end{theorem} 

\subsection{Proofs}\label{subsec:proofs}
\subsubsection{Pre-results}
We would like to introduce first two estimates about the process $X^{t,x}$ defined in \eqref{eq: SDE} which will be used in order to prove the above results. They are related to moments of $X^{t,x}$ under the probabilities $\bP$ and $\bP^{u,v}$, $(u,v)\in \cM$ (see. Karatzas, I.1991 \cite{karatzas}, pp.306). Indeed we have:
\begin{equation}\label{eq: est x}
\forall q\in [1,\infty),\quad \E\big[(\sup_{s\leq T}|X_s^{t,x}|)^{q}\big]\leq C(1+|x|^{q})
\end{equation}
and for any $(u,v)\in \cM$
\begin{equation}\label{eq: est x weak}
\forall q\in [1, \infty),\quad \E^{u,v}_{t,x}\big[(\sup_{s\leq T}|X_s^{t,x}|)^{q}\big]\leq C(1+|x|^{q}).
\end{equation}
Finally let us recall the following important result by U.G.Hausmann (see Theorem 2.2, pp.14 \cite{haussmann}) related to the integrability of the exponential local martingale defined by \eqref{eq:density fun}.
\begin{lemma}\label{lemma: haussmann}({\cite{haussmann}, pp.14}) Let 
$\Theta$ be a $\cP\otimes \cal{B} (\bR)$-measurable application from $[0,T]\times \Omega\times \bR$ to $\bR$ which is of uniformly linear growth, that is, $\bP$-a.s. $\forall(s,x)\in [0,T]\times \bR$, $|\Theta(s,\omega,x)|\leq C_0(1+|x|)$. Then, there exist constants $p\in (1,2)$ and $C$, where $p$ depends only on $C_0$ while the constant $C$, depends only on $p$, but not on $\Theta$, such that:
\[
\E\Big[\Big(\zeta_T\{\Theta(s, X_s^{t,x})\}\Big)^p\Big]\leq C,
\]
where the process $\zeta_T(.)$ is the density function defined by \eqref{eq:density fun}.
\end{lemma}

As a by-product we have:
\begin{corollary} \label{coro: Haussmann}
For any admissible control $(u,v)\in \cM$ and $(t,x)\in [0,T]\times \bR$, there exists a constant $p\in (1,2)$ such that:
\[
\E\Big[\Big(\zeta_T\{\Gamma(s, X_s^{t,x},u_s, v_s)\}\Big)^p\Big]\leq C.
\]
\end{corollary}
\subsubsection{Proof of Proposition \ref{prop: yuv}}\label{subsubsec: proof prop 1}
We will prove this Proposition by constructing the candidate solution of BSDE \eqref{eq:BSDE yiuv} directly. Then we check by It\^o's formula that, the process defined is exactly the solution we anticipate. In this proof, Corollary \ref{coro: Haussmann} plays an important role. Let us illustrate it for player $\pi_1$. The same can be done for 
player $\pi_2$.

For simplicity, only in this proof, we use the notation $(\Y, \Z)$ instead of $(Y^{1;x_0;u,v}, Z^{1;x_0;u,v})$.
\medskip

For any $(u,v)\in \cM$, let us define the process $(\Y_s)_{s\leq T}$ as follows:
  \begin{align}\label{eq:def yuv}
\Y_s:=\E^{u,v}_{x_0}[g_1(X_T^{0,x_0})|\cF_s],\quad \forall s\leq T.
  \end{align}
This process is well defined by noticing that, for any constant $q\geq 1$, we have $\E^{u,v}_{x_0}[|g_1(X_T^{0,x_0})|^{q}]\leq C\E^{u,v}_{x_0}[C(1+\sup_{s\leq T}|X_s^{0,x_0}|^{q\gamma})]<\infty$ which is obtained by (\ref{(A2)}) and \eqref{eq: est x weak}. For writing convenience, we denote by $\zeta_s$, the function\\ $\zeta_s(\Gamma(., X_.^{0,x_0},u_.,v_.))$ as mentioned in \eqref{eq:density fun}. Therefore \eqref{eq:def yuv} can be transformed into:
  \begin{flushleft}
  $\qquad \qquad
  \Y_s=\zeta_s^{-1}\E[\zeta_T\cdot g_1(X_T^{0,x_0})|\cF_s],\quad \forall s\leq T.
  $
  \end{flushleft}
In the following, we show that $\zeta_T\cdot g_1(X_T^{0,x_0})\in L^{\bar{q}}$ for some $\bar{q}\in (1,2)$. Indeed, according to Corollary \ref{coro: Haussmann}, there exists some $p_0\in (1,2)$, such that $\zeta_T\in L^{p_0}(d\bP)$. Therefore, for any $\bar{q}\in (1,p_0)$, Young's inequality leads to:
  \begin{flushleft}
  $\qquad \qquad
  \E[|\zeta_T\cdot g_1(X_T^{0,x_0})|^{\bar{q}}]\leq \frac{\bar{q}}{p_0}\E[|\zeta_T|^{p_0}]+\frac{p_0-\bar{q}}{p_0}\E[|g_1(X_T^{0,x_0})|^{\bar{q}\cdot \frac{p_0}{p_0-\bar{q}}}],
  $
  \end{flushleft}
which is obviously finite by the polynomial growth of $g_1$ and \eqref{eq: est x}. Therefore the process $Y^{u,v}$ is defined. On the other hand, by Doob's inequality, (\ref{(A2)}) and estimate (\ref{eq: est x weak}) we have: 
\be\label{estimyuv}
\forall q>1, \E^{u,v}_{x_0}[\sup_{s\leq T}|Y^{u,v}_s|^q]\leq C\E^{u,v}_{x_0}[|g_1(X_T^{0,x_0})|^q]\leq C(1+|x_0|^q).
\ee

Next thanks to representation Theorem of martingales (\cite{revuzyor}, pp.199) applied to the process $(\E[\zeta_T\cdot g_1(X_T^{0,x_0})|\cF_s])_{s\le T}$, there exists a $\cP$-measurable and $\bR$-valued process $(\Delta_s)_{s\leq T}$ which satisfies
  $
  \E[(\int_0^T|\Delta_s|^2ds)^{\frac{\bar{q}}{2}}]<\infty
  $
and 
  \begin{flushleft}
  $\qquad \qquad
  \Y_s=\zeta_s^{-1}\{\E[\zeta_T\cdot g_1(X_T^{0,x_0})]+\int_0^s \Delta_rdB_r\}:=\zeta_s^{-1}R_s, \quad \forall s\leq T,
  $
  \end{flushleft}
with $R_s:=\E[\zeta_T\cdot g_1(X_T^{0,x_0})|\cF_s]=\E[\zeta_T\cdot g_1(X_T^{0,x_0})]+\int_0^s \Delta_rdB_r$. Next let us set
\begin{equation}\label{eq:def zuv}
  \Z_s:=-\zeta_s^{-1}\big\{R_s\Gamma(s, X_s^{0,x_0},u_s,v_s)-\Delta_s\big\}, \,s\leq T. 
  \end{equation}
As for any $s\in [0,T]$, $$
  d\zeta_s=\zeta_s\cdot \Gamma(s, X_s^{0,x_0},u_s,v_s)dB_s$$ then by It\^o's formula applied to $\zeta^{-1}R$, one can easily check that $(Y^{u,v},Z^{u,v})$ verifies (\ref{eq:BSDE yiuv}) for $i=1$. The same happens for $i=2$.
\ms

Now for any $(u,v)\in \cM$ and $s\leq T$, if $(B^{u,v}_s)_{s\leq T}$ is the Brownian motion under $\bP^{u,v}_{x_0}$, we then deduce from (\ref{eq:BSDE yiuv}) that 
$$
-dY^{i;x_0;u,v}_s=-Z^{i;x_0;u,v}_sdB^{u,v}_s,\,\;Y^{i;x_0;u,v}_T=g_i(X^{0,x_0}_T) 
$$ and $Y^{i;x_0;u,v}_0=\E^{u,v}_{x_0}[g_i(X^{0,x_0}_T]=J_i(u,v)$ since $\cF_0$ is the trivial $\sigma$-algebra completed with the $\bP$-null sets of $\cF$ and taking into account that $\bP$ and $\bP^{u,v}_{x_0}$ are equivalent probabilities. 
Therefore taking into account of (\ref{estimyuv}) and using the Burkholder-Davis-Gundy inequality we have
$$
\forall q>1, \E^{u,v}_{x_0}[(\int_0^T|Z^{i;{x_0};u,v}_r|^2dr)^{\frac{q}{2}}]<\infty.$$
This and (\ref{estimyuv}) imply the estimate (\ref{eq:esti yuv and zuv}) of Proposition \ref{prop: yuv} for $q>1$. Finally for $q=1$, (\ref{eq:esti yuv and zuv}) is obviously true since it is valid for any $q>1$.
The proof of the Proposition \ref{prop: yuv} is completed.\qedh

\subsubsection{Proof of Proposition \ref{prop: equilibrium}}

For $s\leq T$, let us set $\bu_s=\bu(Z_s^1,\theta_s)$ and $\bv_s=\bv(Z_s^2,\vartheta_s)$, then $(\bu,\bv)\in \cM$. On the other hand, thanks to Proposition \ref{prop: yuv}, we obviously have $Y_0^1=J_1(\bu,\bv)$.
\medskip

Next let $u$ be an arbitrary element of $\cM_1$ and let us show that $Y^1\geq Y^{1;x_0;u,\bv}$, which yields $Y^1_0=J_1(\bu,\bv)\geq Y^{1;x_0;u, \bv}_0=J^1(u, \bv)$.
\medskip

The control $(u, \bv)$ is admissible and thanks to Proposition \ref{prop: yuv}, there exists a pair of $\mathcal{P}$-measurable processes $(Y^{1;x_0;u,\bv},Z^{1;x_0;u,\bv})$ such that for any $q>1$,
  \begin{equation}\label{eq:bsde yuvb}
\left\{
   \begin{aligned}
&\E^{u,\bv}_{x_0}\Big[\sup_{0\leq s\leq T}|Y_s^{1;x_0;u,\bv}|^q+\big(\int_0^T|Z_s^{1;x_0;u,\bv}|^2ds\big)^{\frac{q}{2}}\Big]<\infty\,;\\
&Y_s^{1;x_0;u,\bv}=g_1(X_T^{0,x_0})+\int_s^T\! H_1(r, X^{0,x_0}_r, Z_r^{1;x_0;u,\bv}, u_r, \bv_r)dr-\int_s^T\! Z_r^{1;x_0;u,\bv}dB_r,\\
 &\qquad \qquad \qquad \forall s\leq T.
   \end{aligned}
\right.
  \end{equation}

\noindent Afterwards, we aim to compare $Y^1$ in \eqref{eq:main bsde} and $Y^{1;x_0;u,\bv}$ in \eqref{eq:bsde yuvb}. So let us denote by
  \begin{equation*}
\triangle Y=  Y^{1;x_0;u, \bv}-Y^1 \mbox{ and }\triangle Z=  Z^{1;x_0;u,\bv}-Z^1.
  \end{equation*}
\noindent For $k\geq 0$, we define the stopping time $\tau_k$ as follows:
 \begin{center}
 $
\tau_k:= \inf\{s\geq 0, |\triangle
Y_s|+\int_0^s|\triangle Z_r|^2dr\geq k\}\wedge T.
 $
 \end{center}
 
\no The sequence of stopping times $(\tau_k)_{k\geq 0}$ is of stationary type and converges to $T$. Next applying It\^{o}-Meyer formula to $|(\triangle Y)^+|^q$ $(q>1)$ (see Theorem 71, P. Protter, \cite{protter}, pp.221), between $s\wedge \tau_k$ and $\tau_k$, we obtain: $\forall s\leq T$, 
  \begin{equation}\label{eq:tem 3.17}
     \begin{aligned}
|(\triangle &Y_{s\wedge \tau_k})^+|^q+c(q)\int_{s\wedge \tau_k}^{\tau_k}|(\triangle Y_r)^+|^{q-2}1_{\triangle Y_r>0}|\triangle Z_r|^2dr\\
&= |(\triangle Y_{\tau_k})^+|^q+ q\int_{s\wedge
\tau_k}^{\tau_k}|(\triangle Y_{r})^+|^{q-1}
1_{\triangle Y_r>0} \big[H_1(r, X_r^{0,x_0}, Z_r^{1;x_0;u, \bv}, u_r, \bv_r)-\\
&\qquad H_1(r, X_r^{0,x_0}, Z_r^1, \bu_r,
\bv_r)\big]dr -
q\int_{s\wedge \tau_k}^{\tau_k}|(\triangle
Y_{r})^+|^{q-1}1_{\triangle Y_r>0}\triangle Z_rdB_r
  \end{aligned}
  \end{equation}
where $c(q)=\frac{1}{2}q(q-1)$. Besides for any $s\le T$,  
  \begin{align*}
H_1(s, X_s^{0,x_0}, Z_s^{1;x_0;u, \bv}, u_s, \bv_s)&- H_1(s, X_s^{0,x}, Z_s^1, \bu_s,
\bv_s)=\\
&H_1(s, X_s^{0,x_0}, Z_s^{1;x_0;u, \bv}, u_s, \bv_s)- H_1(s,X_s^{0,x_0},Z_s^1, u_s, \bv_s) \\
&+H_1(s, X_s^{0,x_0}, Z_s^1, u_s, \bv_s)-H_1(s, X_s^{0,x_0}, Z_s^1, \bu_s, \bv_s)
  \end{align*}
\no Considering now the generalized Isaacs' condition \eqref{eq:Isaacs}, we have $$H_1(s, X_s^{0,x_0}, Z_s^1, u_s, \bv_s)- H_1(s,X_s^{0,x_0},Z_s^1, \bu_s, \bv_s)\leq 0,\quad \forall s\leq T.$$Additionally
$$H_1(s, X_s^{0,x_0}, Z_s^{1;x_0;u, \bv}, u_s, \bv_s)-H_1(s, X_s^{0,x_0}, Z_s^1, u_s, \bv_s)=\triangle Z_s\Gamma(s,X_s^{0,x_0},u_s,\bv_s).$$Thus equation \eqref{eq:tem 3.17} can be simplified into: $\forall s\in [0,T]$, 
\begin{equation}
\begin{aligned}
|(\triangle &Y_{s\wedge \tau_k})^+|^q+c(q)\int_{s\wedge \tau_k}^{\tau_k}|(\triangle Y_r)^+|^{q-2}1_{\triangle Y_r>0}|\triangle Z_r|^2dr\\
&\leq |(\triangle Y_{\tau_k})^+|^q + q\int_{s\wedge
\tau_k}^{\tau_k}|(\triangle Y_{r})^+|^{q-1}
\mathbbm{1}_{\triangle Y_r>0} \triangle Z_r \Gamma(r, X_r^{0,x_0}, u_r, \bv_r)dr\\
&\quad - q\int_{s\wedge \tau_k}^{\tau_k}|(\triangle Y_{r})^+|^{q-1}1_{\triangle Y_r>0}\triangle Z_rdB_r\nonumber\\
&=|(\triangle Y_{\tau_k})^+|^q-q\int_{s\wedge \tau_k}^{\tau_k}|(\triangle
Y_{r})^+|^{q-1}1_{\triangle Y_r>0}\triangle Z_rdB_r^{u,\bv},
\end{aligned}
  \end{equation}
where $B^{u, \bv}= (B_s- \int_0^{s} \Gamma(r, X_r^{0,x_0}, u_r, \bv_r)dr)_{s\leq T}$ is an $(\cF_s, \bP^{u, \bv}_{x_0})$-Brownian motion. Then for any $s\leq T$,
  \begin{equation*}
\begin{aligned}
|(\triangle Y_{s\wedge \tau_k})^+|^q\leq |(\triangle Y_{\tau_k})^+|^q- q\int_{s\wedge \tau_k}^{\tau_k}|(\triangle Y_{r})^+|^{q-1}1_{\triangle Y_r>0}\triangle Z_rdB_r^{u, \bv}.
\end{aligned}
  \end{equation*}
\noindent By definition of the stopping time $\tau_k$, we have
  \begin{center}
$
\E^{u, \bv}_{x_0}\big[\int_{s\wedge \tau_k}^{\tau_k}|(\triangle Y_r)^+|^{q-1}1_{\triangle Y_r>0}\triangle Z_rdB_r^{u,\bv}\big]=0.
$
\end{center}
Thus for any $s\le T$,
  \begin{align}\label{eq:tem 3.18}
\E^{u, \bv}_{x_0}\left[|(\triangle Y_{s\wedge \tau_k})^+|^q\right] &\leq  \E_{x_0}^{u, \bv}\big[|(Y^{1;x_0;u,\bv}_{\tau_k}-Y_{\tau_k}^1)^+|^q\big].
  \end{align}
Next taking into account \eqref{eq: est x weak} and the fact that $Y^1$ has a representation through $X^{0,x_0}$ and $\eta^1(s,y)$, $(s,y)\in \sp$, a deterministic function with polynomial growth, we have
\begin{equation}\label{eq:uni int}
\E^{u, \bv}_{x_0}\Big[\sup_{s\leq T}(|Y^1_s|+|Y^{1;x_0;u,\bv}_s|)^q\Big]<\infty\end{equation}
As the sequence $((Y_{\tau_k}^{1;x_0;u,\bv}-Y^1_{\tau_k})^+)_{k\geq 0}$
converges to $0$ when $k\rightarrow \infty$,  $\bP_{x_0}^{u,\bv}$-a.s., since $\lim_{k\rightarrow \infty} Y_{\tau_k}^{1;x_0;u,\bv}=\lim_{k\rightarrow \infty}Y_{\tau_k}^1=g_1(X_T^{0,x})\ \bP_{x_0}^{u,\bv}$-a.s.. Then it converges also to $0$ in $L^1(d\bP_{x_0}^{u,\bv})$ thanks to \eqref{eq:uni int}. Take now $k\rightarrow \infty$ on  \eqref{eq:tem 3.18}, it follows from Fatou's Lemma that:
 $$
   \E^{u, \bv}_{x_0}\left[(\triangle Y_{s})^+\right]=0,\,\,\forall s\leq T, 
 $$
which implies that $Y^1\geq Y^{1;x_0;u, \bv}$, $\bP$-a.s., since the probabilities $\bP_{x_0}^{u,\bv}$ and $\bP$ are equivalent. Thus 
 $Y^1_0=J^1(\bu, \bv)\geq Y^{1;x_0;u, \bv}_0=J^1(u, \bv)$. 
\medskip
 
Similarly, we can show that, $Y^2_0=J^2(\bu, \bv)\geq Y^{2;x_0;\bu, v}_0=J^2(\bu,v)$ for arbitrary $v\in \cM_2$. Henceforth $(\bu,\bv)$ is a Nash equilibrium point for the NZSDG. \qedh

\subsubsection{Proof of Theorem \ref{th: main bsde}}\label{subsubsec: proof th}
The proof will be split into several steps. Firstly, we construct an approximating sequence of BSDEs with continuous Lipschitz generators by smoothing the functions $\bar{u}$ and $\bar{v}$. Next we provide appropriate uniform estimates of the solutions of the approximating scheme. Finally we show that the approximating scheme contains at least a convergent subsequence which provides the stochastic processes and deterministic functions verifying the requirements of Proposition \ref{prop: equilibrium}. 
\medskip

\no \underline{\textit{Step 1}}: Approximating scheme.
\ms

At the beginning of this proof, we would like to clarify that the functions $p\in \bR\mapsto p\bar{u}(p,\epsilon_1)$ and $q\in \bR\mapsto q\bar{v}(q,\epsilon_2)$ are uniformly Lipschitz for any $\epsilon_1$ and $\epsilon_2$, since $p\bar{u}(p,\epsilon_1)=p\bar{u}(p,0)=\sup_{u\in U}pu$ and $q\bar{v}(q,\epsilon_2)=q\bar{v}(q,0)=\sup_{v\in V}qv$. Hereafter $\bar{u}(p,0)$ (resp. $\bar{v}(q,0)$) will be simply denoted by $\bar{u}(p)$ (resp. $\bar{v}(q)$).
\medskip

\no Next for integer $n\ge 1$, let $\bar{u}^n$ and $\bar{v}^n$ be the functions defined as follows:
\begin{equation*}
\bar{u}^{n}(p)=\left\{
\begin{aligned}
0&\mbox{ if}\ \ p\leq -1/n,\\
1&\mbox{ if}\ \ p\geq 0,\\
np+1&\mbox{ if}\ \ p\in (-1/n,0), 
\end{aligned}
\right. \mbox{ and }
\quad
\bar{v}^{n}(q)=\left\{
\begin{aligned}
-1&\mbox{ if}\ \ q\leq -1/n,\\
1&\mbox{ if}\ \ q\geq 1/n,\\
nq&\mbox{ if}\ \ q\in (-1/n,1/n).
\end{aligned}
\right.
\end{equation*}
Note that $\bar{u}^n$ and $\bar{v}^n$ are Lipschitz in $p$ and $q$ respectively. Roughly speaking, they are the approximations of $\bar{u}$ and $\bar{v}$. Below,  let $\Phi_n$ be the truncation function $x\in \bR\mapsto \Phi_n(x)=(x\wedge n)\vee(-n) \in \bR$, which is bounded by $n$. Now for $(t,x)\in \sp$ and $n\geq 1$, we consider the following BSDE of dimension two, with Lipschitz generator: For any $s\leq T$, 
  \begin{equation}\label{eq:bsdeapprox}
\left\{
 \begin{aligned}
-dY_s^{1,n;t,x}&= \{\Phi_n(Z^{1,n;t,x}_r)\Phi_n(f(r,X_r^{t,x}))+
\Phi_n(Z^{1,n;t,x}_r\bar{u}(Z^{1,n;t,x}_r))+\\
&\Phi_n(Z^{1,n;t,x}_r)\bar{v}^n(Z^{2,n;t,x}_r)\}dr-Z_r^{1,n;t,x}dB_r ,\  Y_T^{1,n;t,x}=g_1(X_T^{t,x});\\
-dY_s^{2,n;t,x}&=\{\Phi_n(Z^{2,n;t,x}_r)\Phi_n(f(r,X_r^{t,x}))+
\Phi_n(Z^{2,n;t,x}_r\bar{v}(Z^{2,n;t,x}_r))+\\
& \Phi_n(Z^{2,n;t,x}_r)\bar{u}^n(Z^{1,n;t,x}_r)\}dr-Z_r^{2,n;t,x}dB_r,\ Y_T^{2,n;t,x}=g_2(X_T^{t,x}).
 \end{aligned}
\right.
  \end{equation}
From Pardoux-Peng's result (\cite{pardoux peng}), for any $n\ge 1$, this equation has a unique solution $(Y^{i,n;t,x}, Z^{i,n;t,x})\in \mathcal{S}_T^2(\bR)\times \cH_T^2(\bR)$, $i=1,2$. Taking account of the result by El-Karoui et al.(\cite{el karoui}, pp.46, Theorem 4.1), there exist measurable deterministic functions $\eta^{i,n}$ and $\varsigma^{i,n}$, $i=1,2$ and $n\ge 1$, defined on $\sp$ and $\bR$-valued such that:
  \begin{equation}\label{eq:eta in varsigma in}
 \forall s\in [t,T],\,\,Y_s^{i,n;t,x}=\eta^{i,n}(s, X_s^{t,x})\quad \text{and}\quad Z_s^{i,n;t,x}=\varsigma^{i,n}(s, X_s^{t,x}).
  \end{equation}
Moreover, for $n\geq 1$ and $i=1,2,$ the functions $\eta^{i,n}$ verify: 
  \begin{equation}\label{eq:eta n}\forall (t,x) \in [0,T]\times \bR,\,\,
\eta^{i,n}(t,x)=\E[g_i(X_T^{t,x})]+\int_t^T H_i^n(r, X_r^{t,x})dr]
  \end{equation}
with, for any $(s,x)\in [0,T]\times \bR$, 
  \begin{equation}\label{eq:Hn}
\left\{
\begin{aligned}
H_1^n(s,x)&=\Phi_n(\varsigma^{1,n}(s,x))\Phi_n(f(s,x))+
\Phi_n\big(\varsigma^{1,n}(s,x)\bar{u}(\varsigma^{1,n}(s,x))\big)+\\
&\quad + \Phi_n(\varsigma^{1,n}(s,x))\bar{v}^n(\varsigma^{2,n}(s,x));\\
H_2^n(s,x)&=\Phi_n(\varsigma^{2,n}(s,x))\Phi_n(f(s,x))+
\Phi_n\big(\varsigma^{2,n}(s,x)\bar{v}(\varsigma^{2,n}(s,x))\big)+\\
&\quad +\Phi_n(\varsigma^{2,n}(s,x))\bar{u}^n(\varsigma^{1,n}(s,x)).
\end{aligned}
\right.
  \end{equation}
\medskip

\no \underline{\textit{Step 2}}: Estimates for processes $(Y^{i,n;t,x}, Z^{i,n;t,x}),i=1,2$. 
\ms

In order to show the needed uniform estimates for $Y^{i,n;t,x}$ of BSDE \eqref{eq:bsdeapprox}, we use comparison. For that let us consider the following BSDE: For $i=1,2$ and any $s\in[0,T]$, 
\begin{equation}\label{eq: bar yin}
\bar{Y}_s^{i,n}=g_i(X_T^{t,x})+\int_s^T \{\Phi_n(C(1+|X_r^{t,x}|))|\bar{Z}_r^{i,n}|+C|\bar{Z}_r^{i,n}|\}dr-\int_s^T \bar{Z}_r^{i,n}dB_r,
\end{equation}
\no where the constant $C$ is the one such that the generators 
$(H_i^n(s, X_s^{t,x}))_{s\le T}$ satisfy
\be\label{compourestim}\forall s\le T,\,\,|H_i^n(s, X_s^{t,x})|\leq \Phi_n(C(1+|X_s^{t,x}|))|Z_s^{i,n;t,x}|+C|Z_s^{i,n;t,x}|.\ee
This constant exists since $f$ is of linear growth and $\bu$, $\bv$, $\bu^n$ and $\bv^n$ are uniformly bounded. Observing now that the application $z\in \bR\mapsto \Phi_n(C(1+|X_r^{t,x}|))|z|+C|z|$ is Lipschitz continuous, therefore the solution $(\bar{Y}^{i,n},\bar{Z}^{i,n})$ of the above BSDE (\ref{eq: bar yin}) exists in the space $\mathcal{S}_T^2(\bR)\times \cH_T^2(\bR)$ and is unique. Note that we have omitted the dependence w.r.t $(t,x)$ of $(\bar{Y}^{i,n},\bar{Z}^{i,n})$ in order to alleviate notations as there is no possible confusion. On the other hand by the standard comparison theorem of solutions of BSDEs 
(\cite{el karoui}, pp.46, Theorem 4.1) one has
\be\label{comp2}
\bar Y^{i,n}\geq Y^{i,n;t,x}, \bP-a.s.\ee
Next provided that we show uniform estimates, w.r.t. $n$, for $\bar{Y}^{i,n}$, then estimates for $Y^{i,n;t,x}$ will be an immediate consequence. Below, we will focus on the properties of $\bar{Y}^{i,n}$.
\medskip

Using again the result by El Karoui et al.  yields that, there exist deterministic measurable functions $\bar{\eta}^{i,n}:[0,T]\times \bR\rightarrow \bR$ such that, for any $s\in [t,T]$, 
\begin{equation}\label{eq: bar eta in}
\bar{Y}_s^{i,n}=\bar{\eta}^{i,n}(s, X_s^{t,x}),\ i=1,2.
\end{equation}

Next let us consider the process\\ $B^{i,n}=(B_s-\int_0^s[\Phi_n(C(1+|X_r^{t,x}|))+C]\text{sign}(\bar Z_r^{i,n})dr)_{s\leq T}$, $i=1,2$, which is, thanks to Girsanov's Theorem, a Brownian motion under the probability $\bP^{i,n}$ on $(\Omega, \cF)$ whose density with respect to $\bP$ is $\zeta_T\{[\Phi_n(C(1+|X_s^{t,x}|))+C]\text{sign}(\bar Z_s^{i,n})\}$ where for any $z\in \bR, \text{sign}(z)=1_{\{|z|\neq 0\}}\frac{z}{|z|}$ and $\zeta_T(.)$ is defined by \eqref{eq:density fun}. Then the BSDE \eqref{eq: bar yin} will be simplified into,
  \[
\bar{Y}_s^{i,n}=g_i(X_T^{t,x})-\int_s^T \bar{Z}_r^{i,n}dB_r^{i,n},\ s\le T,\ i=1,2.
  \]
\no In view of \eqref{eq: bar eta in}, we obtain,
  \[
\bar{\eta}^{i,n}(t,x)=\E^{i,n}[g_i(X_T^{t,x})|\cF_t],\ i=1,2,
  \]
\no where $\E^{i,n}$ is the expectation under probability $\bP^{i,n}$. By taking the expectation on both sides of the above equation under the probability $\bP^{i,n}$ and considering $\bar{\eta}^{i,n}(t,x)$ is deterministic, we arrive at,
  \[
\bar{\eta}^{i,n}(t,x)=\E^{i,n}[g_i(X_T^{t,x})],\ i=1,2.
  \]
As the functions $g_i$, $i=1,2$, verify the polynomial growth condition (\ref{(A2)}) 
and for any $s\leq T$, 
$$
|\{\Phi_n(C(1+|X_s^{t,x}|))+C\}\text{sign}(\bar Z_s^{i,n})|\le \bar C(1+|X^{t,x}_s|$$
then by estimate \eqref{eq: est x weak} we have, 
for some constant $\lambda \geq 0$ and $C$ a constant which does not depend on $n$, 
$$
|\bar{\eta}^{i,n}(t,x)|\leq C(1+|x|^{\lambda}).
$$ Therefore by (\ref{comp2}) we obtain $\eta^{i,n}(t,x)\le C(1+|x|^\lambda)$, for any $(t,x)\in \sp$. 
In a similar way, we can show that $\eta^{i,n}(t,x)\geq -C(1+|x|^{\lambda})$, $(t,x)\in[0,T]\times \bR$. Therefore, $\eta^{i,n},\ i=1,2$ are of polynomial growth with respect to $(t,x)$ uniformly in $n$.

To conclude this step, we have the following results: There exists a constant $C$ independent of $n$ and $t,x$ such that, for $(t,x)\in [0,T]\times \bR$, $i=1,2,$
   \begin{equation}\label{eq:rst step 2}
   \left\{
     \begin{aligned}
     &\text{(a) } |\eta^{i,n}(t,x)|\leq C(1+|x|^{\lambda}),\text{ for any }\lambda\geq 0;\\
     &\text{(b) by the combination of (a), \eqref{eq: est x} with \eqref{eq:eta in varsigma in}, it holds:  }\\&\qquad  \forall \alpha\geq 1,\,\,
     \E[\sup_{s\in[t,T]}|Y^{i,n;t,x}_s|^{\alpha}]\leq C\,;\\
     &\text{(c) } \text{for any } (t,x)\in \sp,\  \E[\int_t^T |Z_s^{i,n;t,x}|^2ds]\leq C \text{ which is a }\\
     &\qquad \text{straightforward result by using It\^o's formula with } (Y^{i,n;t,x})^2 \mbox{ and using (b)}.
     \end{aligned}
   \right.
   \end{equation}
\no \underline{\textit{Step 3}}: Convergence of a subsequence of $(Y^{i,n;0,x},Z^{i,n;0,x})_{n\geq 1}$, $i=1,2$.
\ms

First let us define on $\bR$ the measure $\mu(0,x_0;s,dy)$ as the law under $\bP$ of $X_s^{0,x_0}$, i.e., $\mu(0,x_0;s,dy):=\bP(X^{0,x_0}_s\in dy)
=\frac{1}{\sqrt{2\pi s}}e^{-\frac{(y-x_0)^2}{2s}}dy$ for any $s\in (0,T]$.

Let $q\in (1,2)$ be fixed. We are going to show that the sequence $(H_i^n(s, y))_{n\ge 1}$ belongs to $L^q([0,T]\times \bR; \mu(0,x_0;s,dy)ds)$, $i=1,2$. Actually,
  \begin{align}\label{eq: Hin in lq}
\E[\int_0^T|H_i^n(s,X_s^{0,x_0})|^q ds]
&=\int_{[0,T]\times \bR}|H_i^n(s,y)|^q\mu(0,x_0;s,dy)ds\nonumber\\
&\leq C\E[\int_0^T|Z^{i,n;0,x_0}_s|^q(1+|X_s^{0,x}|^q)ds]\nonumber\\
&\leq C\{\E[\int_0^T|Z_s^{i,n;0,x_0}|^2ds]+\E[1+\sup_{s\in[0,T]}|X_s^{0,x_0}|^{\frac{2q}{2-q}}]\}\le C\nonumber\\
\end{align}
($C$ is a generic constant whose value may change from line to line). The last inequality is obtained from the fact that $\E[\int_0^T |Z_s^{i,n;0,x_0}|^2 ds]\leq C$ and estimate \eqref{eq: est x}. As a result, there exists a subsequence $\{n_k\}$ (still denoted by $\{n\}$ for simplification) and two $\mathcal{B}([0,T]\times \bR)$-measurable deterministic functions $H_i(s,y)$, $i=1,2$, such that, 
  \begin{equation}\label{eq:Hin weak conv}
H_i^n\rightharpoonup H_i \text{ weakly in } L^q([0,T]\times \bR; \mu(0,x_0;s,dy)ds).
  \end{equation}
Next we focus on passing from the weak convergence to strong sense convergence by proving that $(\eta^{i,n}(t,x))_{n\geq 1}$ defined in \eqref{eq:eta n} is a Cauchy sequence for each $(t,x)\in [0,T]\times \bR$, $i=1,2.$ Let $(t,x)$ be fixed in $[0,T)\times \bR$ (w.l.o.g we assume $t<T$), $\delta>0, k,n$ and $m\geq 1$ be integers. From \eqref{eq:eta n}, we have,
  \begin{align}\label{eq:eta in - eta im}
|\eta^{i,n}(t,x)-\eta^{i,m}(t,x)|&= \big|\E[\int_t^T H_i^n(s, X_s^{t,x})-H_i^m(s,X_s^{t,x})ds]\big|\nonumber\\
&\leq \big|\E[\int_t^{t+\delta} H_i^n(s, X_s^{t,x})-H_i^m(s,X_s^{t,x})ds]\big|\nonumber\\
&\quad +\big|\E[\int_{t+\delta}^T (H_i^n(s, X_s^{t,x})-H_i^m(s,X_s^{t,x}))\cdot 1_{\{|X_s^{t,x}|\leq k\}}ds]\big|\nonumber\\
&\quad +\big|\E[\int_{t+\delta}^T (H_i^n(s, X_s^{t,x})-H_i^m(s,X_s^{t,x}))\cdot 1_{\{|X_s^{t,x}|> k\}}ds]\big|.\nonumber\\
&
  \end{align}
We first deal with the first term of the right-hand side of \eqref{eq:eta in - eta im}. By H\"older's inequality, definition of $H_i^n(s, X_s^{t,x})$ and (\ref{eq:rst step 2}) we have$$\begin{array}{ll}
    \big|\E[\int_t^{t+\delta}\!\!\! H_i^n(s, X_s^{t,x})-\!\!\!{}&\!\!\!H_i^m(s,X_s^{t,x})ds]\big|\\ &\leq\!   
    \E[\int_t^{t+\delta}| H_i^n(s, X_s^{t,x})-H_i^m(s,X_s^{t,x})|ds] \\{}&\leq
    C\E[\int_t^{t+\delta}|Z^{i,n;t,x}_s|(1+|X^{t,x}_s|)ds]\\{}&\leq  C\delta^{\frac{1}{4}}\{\E[\int_t^T\!\! |Z^{i,n;t,x}_s|^2ds]\}^{\frac{1}{2}}\{\E[\int_t^T\!\! (1+|X^{t,x}_s|)^4ds]\}^{\frac{1}{4}} \\{}& \leq C\delta^{\frac{1}{4}}.
    \end{array}
    $$
We now focus on the third term of the right hand-side of inequality \eqref{eq:eta in - eta im}. In the same way as previous we have: 
$$\begin{array}{ll}
    \big|\E&[\int_{t+\delta}^T (H_i^n(s, X_s^{t,x})-H_i^m(s,X_s^{t,x}))\cdot 1_{\{|X_s^{t,x}|> k\}}ds]\big|\\ &\leq \E[\int_{t+\delta}^T |H_i^n(s, X_s^{t,x})-H_i^m(s,X_s^{t,x})|\cdot 1_{\{|X_s^{t,x}|> k\}}ds]\\     
    &\leq C\{\E[\int_{t+\delta}^T1_{\{|X_s^{t,x}|>k\}}ds]\}^{\frac{1}{4}}
\{\E[\int_t^T\!\! |Z^{i,n;t,x}_s|^2ds]\}^{\frac{1}{2}}\{\E[\int_t^T\!\! (1+|X^{t,x}_s|)^4ds]\}^{\frac{1}{4}} \\     
    &\leq Ck^{-\frac{1}{4}}
    \end{array}
    $$
    by using Markov inequality. Finally we deal with the second term of the right-hand side of \eqref{eq:eta in - eta im}. By using the law of $X^{t,x}_s$ we have
    \begin{align}\label{eq:sec term}
   & \E[\int_{t+\delta}^T (H_i^n(s, X_s^{t,x})-H_i^m(s,X_s^{t,x}))\cdot 1_{\{|X_s^{t,x}|\leq k\}}ds]\nonumber\\&=   
     \int_{\bR}\int_{t+\delta}^T (H_i^n(s,y)-H_i^m(s,y))1_{\{|y|\leq k\}}\frac{1}{\sqrt{2\pi(s-t)}}e^{-\frac{(y-x)^2}{2(s-t)}}dsdy\nonumber\\
    &=\int_{\bR}\int_{t+\delta}^T (H_i^n(s,y)-H_i^m(s,y))1_{\{|y|\leq k\}}\Phi_{t,x,x_0}(s,y)\frac{1}{\sqrt{2\pi s}}e^{-\frac{(y-x_0)^2}{2s}}dy
    \end{align}
where 
$$
\Phi_{t,x,x_0}(s,y)=\frac{\sqrt{s}}{\sqrt{s-t}}e^{-\frac{(y-x)^2t}{2(s-t)}+\frac{(y-x_0)^2}{2s}}
$$which is bounded when $s\in [t+\delta,T]$ and $|y|\leq k$ and then \\
$\Phi_{t,x,x_0}(s,y)1_{[t+\delta,T]\times [-k,k]}(s,y)$ belongs 
$L^{\bar q}([0,T]\times \bR; \mu(0,x_0;s,dy)ds)$ with $\bar q$ is the conjugate of $q$. Now as the sequence of functions $(H_i^n)_{n\geq 1}$ converges weakly in $
L^q([0,T]\times \bR; \mu(0,x_0;s,dy)ds)$ then $$\begin{array}{ll}
 &\E[\int_{t+\delta}^T (H_i^n(s, X_s^{t,x})-H_i^m(s,X_s^{t,x}))\cdot 1_{\{|X_s^{t,x}|\leq k\}}ds]\\&\quad =
 \int_{\bR}\int_{t+\delta}^T (H_i^n(s,y)-H_i^m(s,y))1_{\{|y|\leq k\}}\Phi_{t,x,x_0}(s,y)\frac{1}{\sqrt{2\pi s}}e^{-\frac{(y-x_0)^2}{2s}}dy\\&\quad \rightarrow 0 \mbox{ as }n,m\rightarrow \infty.\end{array}
$$
Henceforth for $i=1,2$, and any $(t,x)\in \sp$, the sequences $(\eta^{i,n}(t,x))_{n\geq 1}$ is of Cauchy type and then converges to a deterministic measurable function $\eta^i(t,x)$. Additionally from the uniform polynomial growth of $\eta^{i,n}$ (see (a) of Step 2) we deduce that $\eta^i$ is also of polynomial growth, i.e., \be\label{polygrowthetai}
\forall (t,x)\in \sp, \,\,|\eta^i(t,x)|\le C(1+|x|^\lambda) \,\,(\lg \ge 0).\ee
Next it turns out that, for $i=1,2$ and any $s\in [0,T]$,
    $$
    \lim_{n\rightarrow \infty} Y_s^{i,n;0,x_0}(\omega)=\eta^i(s,X_s^{0,x_0}(\omega)) \mbox{ and }|Y_s^{i,n;0,x_0}(\omega)|\leq C(1+|X_s^{0,x_0}(\omega)|^{\lambda}),\ \bP-a.s.
    $$
Setting $Y^i:=(\eta^i(s,X_s^{0,x_0}))_{s\leq T}$, then by Lebesgue's dominated convergence theorem it holds \begin{align}\label{eq:yin conv l alpha}
    \E[\int_0^T |Y_s^{i,n;0,x_0}-Y_s^i|^{\alpha}ds]\rightarrow 0, \quad \text{as}\ n     \rightarrow \infty \quad \text{for any}\ \alpha\geq 1\,\, (i=1,2).
    \end{align}

It remains to show the convergence of sequences $((Z_s^{i,n;0,x_0})_{s\leq T})_{n\geq 1}$, $i=1,2$. Taking It\^o's formula with $(Y^{i,n;0,x_0}-Y^{i,m;0,x_0})^2$ and considering (\ref{lg}), we get:
  \begin{align}\label{eq:2.31}
&|Y_s^{i,n;0,x_0}-Y_s^{i,m;0,x_0}|^2+\int_s^T|Z_r^{i,n;0,x_0}-Z_r^{i,m;0,x_0}|^2ds\nonumber\\
&\leq 2\int_s^TC|Y_r^{i,n;0,x_0}-Y_r^{i,m;0,x_0}|(1+|X_r^{0,x_0}|)(|Z_r^{i,n;0,x_0}|+|Z_r^{i,m;0,x_0}|)dr\nonumber\\
&\quad -2\int_s^T(Y_r^{i,n;0,x_0}-Y_r^{i,m;0,x_0})(Z_r^{i,n;0,x_0}-Z_r^{i,m;0,x_0})dB_r,\,\,\forall s\le T.
  \end{align}
\noindent Since for any $a,b,c\in \bR$ and for any $\epsilon>0$, $|abc|\leq \frac{\epsilon^2}{2}a^2+\frac{\epsilon^4}{4}b^4+\frac{1}{4\epsilon^8}c^4$, we then have: $\forall s\le T$, 
  \begin{align}\label{eq:2.32}
  &|Y_s^{i,n;0,x_0}-Y_s^{i,m;0,x_0}|^2+\int_s^T|Z_r^{i,n;0,x_0}-Z_r^{i,m;0,x_0}|^2dr\nonumber\\
  &\leq C\big\{\frac{\epsilon^2}{2}\int_s^T(|Z_r^{i,n;0,x_0}|+|Z_r^{i,m;0,x_0}|)^2dr+\frac{\epsilon^4}{4}\int_s^T(1+|X_r^{0,x_0}|)^4dr\nonumber\\
  &\qquad+\frac{1}{4\epsilon^8}\int_s^T|Y_r^{i,n;0,x_0}-Y_r^{i,m;0,x_0}|^4dr\big\}\nonumber\\
  &\quad-2\int_s^T(Y_r^{i,n;0,x_0}-Y_r^{i,m;0,x_0})(Z_r^{i,n;0,x_0}-Z_r^{i,m;0,x_0})dB_r.
  \end{align}
Since $\epsilon$ is arbitrary, taking now $s=0$, expectation on both hand-sides and the limit w.r.t. $n,m$, combining with \eqref{eq:yin conv l alpha}, \eqref{eq: est x}, \eqref{eq:rst step 2}-(c) yields that, 
   \begin{equation}\label{eq:2.33}
   \limsup_{n,m\rightarrow \infty} \E[\int_0^T |Z_r^{i,n;0,x_0}-Z_r^{i,n;0,x_0}|^2dr]\rightarrow 0,\quad i=1,2.
   \end{equation}
Consequently, for $i=1,2$, the sequence $(Z^{i,n;0,x_0}=(\varsigma^{i,n}(t, X_t^{0,x}))_{t\leq T})_{n\geq 1}$ is convergent in $\cH^2_T(\bR)$ to a process $Z^i$ which belongs also to $\cH^2_T(\bR)$.
On the other hand one can substract a subsequence which we still denote by $\{n\}$ such that 
$(Z^{i,n;0,x_0}_s)_{n\ge 1}\rw Z^i_s$, $ds\otimes d\bP$-a.e. and $\sup_{n\ge 1}|Z^{i,n;0,x_0}_s|$ belongs to $L^2(\sp,ds\otimes d\bP)$. 

Next going back to inequality \eqref{eq:2.32}, taking the supremum on interval $[0,T]$ and using BDG's inequality, we deduce that,
  \begin{align*}
  &\E\big[\sup_{s\in [0,T]}|Y_s^{i,n;0,x_0}-Y_s^{i,m;0,x_0}|^2+\int_0^T|Z_r^{i,n;0,x_0}-Z_r^{i,m;0,x_0}|^2dr\big]\nonumber\\
  &\leq C\E\big\{\frac{\epsilon^2}{2}\int_0^T(|Z_r^{i,n;0,x_0}|+|Z_r^{i,m;0,x_0}|)^2dr+\frac{\epsilon^4}{4}\int_0^T(1+|X_r^{0,x_0}|)^4dr\nonumber\\
  &\qquad \qquad +\frac{1}{4\epsilon^8}\int_0^T|Y_r^{i,n;0,x_0}-Y_r^{i,m;0,x_0}|^4dr\big\}\nonumber\\
  &\quad+\frac{1}{2}\E\big[\sup_{r\in[0,T]}|Y_r^{i,n;0,x_0}-Y_r^{i,m;0,x_0}|^2\big]+2\E\big[\int_0^T|Z_r^{i,n;0,x_0}-Z_r^{i,m;0,x_0}|^2dr\big],
  \end{align*}
\noindent which implies,
  $$
\limsup_{n,m\rightarrow \infty}\E\big[\sup_{s\in[0,T]}|Y_s^{i,n;0,x_0}-Y_s^{i,m;0,x_0}|^2\big]=0,
  $$
since $\epsilon$ is arbitrary and the facts of \eqref{eq: est x}, \eqref{eq:2.33}, \eqref{eq:yin conv l alpha} and \eqref{eq:rst step 2}-(c). Thus the sequence of processes $(Y^{i,n;0,x_0})_{n\geq 1}$ converges in $\mathcal{S}_T^2(\bR)$ to $Y^i$ for $i=1,2$ which are continuous processes.
\medskip

To summarize this step, we have the following results (at least for a subsequence $\{n\}$): for $i=1,2$,
  \begin{equation}\label{eq:rst step3}
\left\{
    \begin{aligned}
    &\text{(a)}\ H_i^n(s,y)\in L^q([0,T]\times \bR;\ \mu(0,x_0;s,dy)ds)\text{ uniformly w.r.t. }n;\\
    &\text{(b)}\ Y^{i,n;0,x_0}\rightarrow_{n\rightarrow \infty} Y^i \text{ in } L^{\alpha}([0,T]\times \bR, ds\otimes d\bP)\text{ for any } \alpha\geq 1,\text{ besides, }\\
    &\qquad Y^{i,n;0,x_0}\rightarrow_{n\rightarrow \infty} Y^i \text{ in } \mathcal{S}_T^2(\bR);\\
    &(c)\  Z^{i,n;0,x_0}\rightarrow_{n\rightarrow\infty}Z^i\text{ in } L^2([0,T]\times \bR,ds\otimes d\bP), \text{ additionally, there exists a}\\
    &\qquad \text{subsequence } \{n\} \text{ s.t }Z^{i,n;0,x_0}\rightarrow_{n\rightarrow\infty}Z^i\ ds\otimes d\bP-a.e.\text{ and }\\
    &\qquad \sup_{n\geq 1}|Z^{i,n;0,x_0}|\in L^2([0,T]\times \bR,ds\otimes d\bP).
    \end{aligned}
\right.
  \end{equation}

\no \underline{\textit{Step 4}}: Convergence of $(H_i^n)_{n\geq 1},\ i=1,2$.
\medskip

In this step, we are going to define the processes $(\theta_s)_{s\leq T}$ and $(\vartheta_s)_{s\leq T}$ and verify that $(Y^i,Z^i)$, $i=1,2$, and $\theta$, $\vartheta$ satisfy (\ref{eq:main bsde}). We  demonstrate first for $i=1$. Let us consider the subsequence which satisfies (\ref{eq:rst step3}).
\medskip

%
Recall \eqref{eq:Hn} for $(t,x)=(0,x_0)$ which reads as: For any $\sT$,
 \begin{equation}\label{eq:rec Hn}
 \begin{aligned}
H_1^n(s,X^{0,x_0}_s)&=\Phi_n(Z_s^{1,n;0,x_0})\Phi_n(f(s,X^{0,x_0}_s))+
\Phi_n(Z_s^{1,n;0,x_0}\bar{u}(Z_s^{1,n;0,x_0}))\nonumber\\
&\quad+\Phi_n(Z_s^{1,n;0,x_0})\bar{v}^n(Z_s^{2,n;0,x_0}). 
\end{aligned}
\end{equation}
Note that, 
 \begin{align*}
 \Phi_n(Z_s^{1,n;0,x_0})\Phi_n(f(s,X^{0,x_0}_s))&+
\Phi_n(Z_s^{1,n;0,x_0}\bar{u}(Z_s^{1,n;0,x_0}))\\
&\rightarrow_{n\rightarrow\infty} Z_s^{1}f(s,X^{0,x_0}_s)+Z_s^1\bar{u}(Z_s^{1}),\  
 ds\otimes d\bP\text{-a.e.}
\end{align*}
since $Z^{1,n;0,x_0}\rightarrow_{n\rightarrow \infty}Z^1$, $ds\otimes d\bP$-a.e. as stated in \eqref{eq:rst step3}-(c), $\Phi_n(x)\rightarrow_{n\rightarrow \infty}x$ and finally by the continuity of $p\in \bR\mapsto p\bar{u}(p)$. This convergence holds also in $\cH^1_T(\bR)$, by Lebesgue's dominated convergence theorem, since 
the process $(\sup_{n\ge 1}|Z_s^{1,n;0,x_0}|)_{\sT}$ belongs to $\cH^2_T(\bR)$, $f$ is of linear growth and $\bu$ is uniformly bounded. The rest part in \eqref{eq:rec Hn} is
 $$\begin{array}{l}
\Phi_n(Z_s^{1,n;0,x_0})\bar{v}^n(Z_s^{2,n;0,x_0})\\\qquad =\Phi_n(Z_s^{1,n;0,x_0})\bar{v}^n(Z_s^{2,n;0,x_0})1_{\{Z_s^{2}\neq 0\}}+\Phi_n(Z_s^{1,n;0,x_0})
\bar{v}^n(Z_s^{2,n;0,x_0})1_{\{Z_s^{2}=0\}}.
\end{array} $$
But $$
 \Phi_n(Z_s^{1,n;0,x_0})\bar{v}^n(Z_s^{2,n;0,x_0})
 1_{\{Z_s^{2}\neq 0\}} \rightarrow_{n\rightarrow \infty} Z_s^1\bar{v}(Z_s^2)1_{\{Z_s^2\neq 0\}} \quad ds\otimes d\bP\text{-a.e.}
 $$
since for any $z\neq 0$, if $n$ is large enough then $\bar v^n(z')=\bar v(z')$ if $z'\in (z-a,z+a)\subset (-\infty,0)\cup (0,+\infty)$ for small $a>0$ and $\bar v$ is continuous in $z$. Once more 
the convergence holds in $\cH^2_T(\bR)$, by Lebesgue's dominated convergence theorem, since 
the process $(\sup_{n\ge 1}|Z_s^{1,n;0,x_0}|)_{\sT}$ belongs to $\cH^2_T(\bR)$ and $\bv^n$ is uniformly bounded. 

To proceed let us define a $\cP$-measurable process $(\vartheta_s)_{s\leq T}$ valued on $V$ as the weak limit in $\cH_T^2(\bR)$ of some subsequence $(\bar{v}^{n_k}(Z^{2,n_k;0,x_0})1_{\{Z^{2}=0\}})_{k\geq 0}$. The weak limit exists since $(\bar{v}^{n_k})_{k\geq 0}$ is bounded. Let now $\tau$ be an arbitrary stopping time such that $\tau \in [0,T]$, $ \bP-a.s.$, then
\begin{align*}
 \int_0^{\tau}&\Phi_{n_k}(Z_s^{1,n_k;0,x_0})\bar{v}^{n_k}(Z_s^{2,n_k;0,x_0})1_{\{Z_s^2=0\}}ds\\
 &\rightarrow_{k\rightarrow\infty} \int_0^{\tau}Z_s^1\vartheta_s1_{\{Z_s^2=0\}}ds\text{ weakly in } L^2(\Omega, d\bP).\end{align*}
Indeed 
\begin{align*}
  \int_0^{\tau}&\Phi_{n_k}(Z_s^{1,n_k;0,x_0})\bar{v}^{n_k}(Z_s^{2,n_k;0,x_0})1_{\{Z_s^2=0\}}ds\\
  &=\quad\int_0^{\tau}(\Phi_{n_k}(Z_s^{1,n_k;0,x_0})-Z_s^1)\bar{v}^{n_k}(Z_s^{2,n_k;0,x_0})1_{\{Z_s^2=0\}}ds\\
  &\quad+\int_0^{\tau}Z_s^1\bar{v}^{n_k}(Z_s^{2,n_k;0,x_0})1_{\{Z_s^2=0\}}ds.
  \end{align*}
\no On the right side, the first integral converges to $0$ in $L^2(d\bP)$ by Lebesgue's dominated convergence theorem since
$\Phi_{n_k}(Z^{1,n_k;0,x_0})\rightarrow Z^1 \ dt\otimes d\bP$-a.e., $\sup_{k\geq 0}|Z_t^{1,n_k}|\in L^2([0,T]\times \bR)$ as shown in \eqref{eq:rst step3}-(c), $Z^1\in L^2([0,T]\times \bR)$, the sequence $(\bv^{n_k})_{k\geq 0}$ is bounded and $|\Phi_{n_k}(x)|\leq |x|$, $\forall x\in \bR$. Below, we will give the weak convergence in $L^2(\Omega, d\bP)$ of the integral $\int_0^{\tau}Z_s^1\bar{v}^{n_k}(Z_s^{2,n_k;0,x_0})1_{\{Z_s^2=0\}}ds$ to $\int_0^{\tau}Z_s^1\vartheta_s1_{\{Z_s^2=0\}}ds$. That is, for any random variable $\xi\in L^2(\Omega, \cF_T, d\bP)$, we need to show, 
  \begin{equation}\label{eq:weak cov int xi}
\E[\xi\int_0^{\tau}Z_s^1\bar{v}^{n_k}(Z_s^{2,n_k;0,x_0})1_{\{Z_s^2=0\}}ds]\rightarrow_{k\rightarrow\infty}\E[\xi\int_0^{\tau}Z_s^1\vartheta_s 1_{\{Z_s^2=0\}}ds]. 
  \end{equation}
\noindent Thanks to martingale representation theorem, there exists a process $( \Lambda_s)_{s\leq T}\in \cH_{T}^2(\bR)$ such that, 
$
\E[\xi|\cF_{\tau}]=\E[\xi]+\int_0^{\tau}\Lambda_sdB_s.
$
Therefore, 
  \begin{align*}
\E\Big[\xi\int_0^{\tau}&Z_s^1\bar{v}^{n_k}(Z_s^{2,n_k;0,x_0})1_{\{Z_s^2=0\}}ds\Big]\\
&=\E\Big[\E\big[\xi\int_0^{\tau}Z_s^1\bar{v}^{n_k}(Z_s^{2,n_k;0,x_0})1_{\{Z_s^2=0\}}ds|\cF_{\tau}]\Big]\\
&=\E\Big[\E[\xi|\cF_{\tau}]\cdot \int_0^{\tau}Z_s^1\bar{v}^{n_k}(Z_s^{2,n_k;0,x_0})1_{\{Z_s^2=0\}}ds\Big]\\
&=\E\Big[\E[\xi]\int_0^{\tau}Z_s^1\bar{v}^{n_k}(Z_s^{2,n_k;0,x_0})1_{\{Z_s^2=0\}}ds\Big]\\
&\qquad+ \E\Big[\int_0^{\tau}\Lambda_sdB_s\cdot\int_0^{\tau}Z_s^1\bar{v}^{n_k}(Z_s^{2,n_k;0,x_0})1_{\{Z_s^2=0\}}ds\Big].
  \end{align*}
Notice that $\E[\xi]\E[\int_0^{\tau}Z_s^1\bar{v}^{n_k}(Z_s^{2,n_k;0,x_0})1_{\{Z_s^2=0\}}ds] \rightarrow_{k\rightarrow \infty} \E[\xi]\E[\int_0^{\tau}Z_s^1\vartheta_s1_{\{Z_s^2=0\}}ds]$, since $(Z_s^1)_{s\leq T}\in \mathcal{H}_T^2(\bR)$ and $\bar{v}^{n_k}(Z_s^{2,n_k;0,x_0})1_{\{Z_s^2=0\}}\rightarrow_{k\rightarrow\infty}\vartheta_s$ weakly in $\mathcal{H}_T^2(\bR)$. Next, by It\^o's formula,
  \begin{align*}
\E[&\int_0^{\tau}\Lambda_sdB_s\cdot\int_0^{\tau}Z_s^1\bar{v}^{n_k}(Z_s^{2,n_k;0,x_0})1_{\{Z_s^2=0\}}ds]\\
&=\E[\int_0^{\tau}\big(\int_0^s\Lambda_udB_u\big)Z_s^1\bar{v}^{n_k}(Z_s^{2,n_k;0,x_0})1_{\{Z_s^2=0\}}ds]+\\
&\quad+\E[\int_0^{\tau}\big(\int_0^sZ_u^1\bar{v}^{n_k}(Z_u^{2,n_k;0,x_0})1_{\{Z_u^2=0\}}du\big)\Lambda_sdB_s].
  \end{align*}
The latter one on the right side is 0, since $\int_0^{.}(\int_0^sZ_u^1\bar{v}^{n_k}(Z_u^{2,n_k;0,x_0})1_{\{Z_u^2=0\}}du)\Lambda_sdB_s$ is an $\cF_{t}$-martingale which is followed by $(Z_s^1)_{s\leq T}\in \cH_T^2(\bR),\ (\Lambda_s)_{s\leq T}\in \cH_T^2(\bR)$ and the boundness of $\bv^{n_k}$ and then \\
$\E[\{\int_0^T|\int_0^sZ_u^1\bar{v}^{n_k}(Z_u^{2,n_k;0,x_0})1_{\{Z_u^2=0\}}du|^2|\Lambda_s|^2ds\}^{\frac{1}{2}}]<\infty$. For the former part, let us denote $\int_0^s\Lambda_udB_u$ by $\psi_s$ for any $s\in [0,\tau]$. Then for any integer $\kappa>0$, we have,
  \begin{align*}
|\E&[\int_0^{\tau}\psi_sZ_s^1\Big(\bar{v}^{n_k}(Z_s^{2,n_k})1_{\{Z_s^2=0\}}-\vartheta_s\Big)ds]|\\
&=|\E[\int_0^{\tau}\psi_sZ_s^1\Big(\bar{v}^{n_k}(Z_s^{2,n_k;0,x_0})1_{\{Z_s^2=0\}}-\vartheta_s\Big)1_{\{|\psi_sZ_s^1|\leq \kappa\}}ds]|+\\%
&\quad+|\E[\int_0^{\tau}\psi_sZ_s^1\Big(\bar{v}^{n_k}(Z_s^{2,n_k;0,x_0})1_{\{Z_s^2=0\}}-\vartheta_s\Big)1_{\{|\psi_sZ_s^1|> \kappa\}}ds]|.
  \end{align*}
On the right side of the above equation, the first component converges to 0 which is the consequence of $\bar{v}^{n_k}(Z_s^{2,n_k;0,x_0})1_{\{Z_s^2=0\}}\rightarrow_{k\rightarrow\infty}\vartheta_s$ weakly in $\cH_T^2(\bR)$. For the second term, considering both $(\bar{v}^{n_k}(Z_s^{2,n_k;0,x_0}))_{s\leq \tau}$ and $(\vartheta_s)_{s\leq \tau}$ are bounded, it is smaller than $ C|\E[\int_0^{\tau}|\psi_sZ_s^1|1_{\{|\psi_sZ_s^1|\geq \kappa\}}ds]|$ which obviously converges to 0 as $\kappa\rightarrow \infty$. Thus \eqref{eq:weak cov int xi} holds true. 

Finally, we also have
 \begin{equation}\label{eq:conv z tau} 
 \int_0^{\tau}Z_s^{1,n_k;0,x_0}dB_s \rightarrow_{k\rightarrow\infty} \int_0^{\tau}Z_s^1dB_s \text{ in } L^2(\Omega, d\bP),
 \end{equation}
 which is obtained from the convergence of $(Z^{1,n_k;0,x_0})_{k\geq 0}$ to $Z^1$ in $\cH_T^2(\bR)$.
 
Then by observing the approximation BSDE \eqref{eq:bsdeapprox} in a forward way, \textit{i.e.} for any stopping time $\tau$,
  $$
Y_{\tau}^{1,n_k;0,x_0}=Y_0^{1,n_k;0,x_0}-\int_0^{\tau} H_1^{n_k}(s, X_s^{0,x_0})ds+\int_0^{\tau} Z_s^{1,n_k;0,x_0}dB_s,
  $$
\no  combining with the convergence of $(Y^{1,n_k;0,x_0})_{k\geq 0}$ to $Y^1$ in $\mathcal{S}_T^2(\bR)$ we infer that
  \begin{equation*}
\bP\text{-a.s.},\ Y_{\tau}^1=Y_0^1-\int_0^{\tau}H_1^*(s,X^{0,x_0}_s,Z_s^1,Z_s^2,\vartheta_s)ds+\int_0^{\tau}Z_s^1dB_s
  \end{equation*}
since 
$$
\int_0^{\tau} H_1^{n_k}(s, X_s^{0,x_0})ds \rightharpoonup_{k\rw \infty}\int_0^{\tau}H_1^*(s,X^{0,x_0}_s,Z_s^1,Z_s^2,\vartheta_s)ds \mbox{ weakly in }L^2(\Omega, d\bP).
$$
As $\tau$ is arbitrary then the processes $Y^1$ and $Y_0^1-\int_0^\cdot H_1^*(s,X^{0,x_0}_s,Z_s^1,Z_s^2,\vartheta_s)ds+\int_0^\cdot Z_s^1dB_s$ are indistinguishable, \textit{i.e.}, $\bP$-a.s.
  \begin{equation*}
\forall s\leq T,\quad Y_s^1=Y_0^1-\int_0^sH_1^*(r,X^{0,x_0}_r,Z_r^1,Z_r^2,\vartheta_r)dr+\int_0^sZ_r^1dB_r.
  \end{equation*}
On the other hand, $Y_T^1=g_1(X^{0,x_0}_T)$, then,
  \begin{equation*}\bP\text{-a.s.},\,\quad
\forall s\leq T,\quad Y_s^1=g_1(X^{0,x_0}_T)+\int_s^TH_1^*(r,X^{0,x_0}_r,Z_r^1,Z_r^2,\vartheta_r)dr-\int_s^TZ_r^1dB_r.
  \end{equation*}
Similarly, for player $\pi_2$, there exists a $\cP$-measurable process $(\theta_s)_{s\leq T}$ valued on $U$, which is obtained in the same way as previously as a weak limit of 
a subsequence of 
$({\bar u}^n(Z^{1,n;0,x_0}_s)1_{\{Z^1_s\neq 0\}})_{s\leq T} \mbox{ in }\cH^2_T(\bR)$ such that 
  \begin{equation*}
\bP\text{-a.s.},\quad \forall s\leq T,\quad Y_s^2=g_2(X^{0,x_0}_T)+\int_s^T
H_2^*(r,X^{0,x_0}_r,Z_r^1,Z_r^2,\theta_r)dr-\int_s^TZ_r^2dB_r.
  \end{equation*}
The proof is completed.\qedh
\section{Generalizations}\label{sec:gener}
In this Section, we are going to deal with some generalizations of Theorem \ref{th: nash} in the following three aspects:
\medskip

\noindent \textbf{(i)} For the drift term $\Gamma$ in SDE \eqref{eq: SDE with Gamma} which reads,
  $$
  \Gamma(t,x,u,v)=f(t,x)+u+v, \,\,u\in [0,1]\mbox{ and }v\in [-1,1],
  $$
one can replace:

(a) $[0,1]$ and $[-1,1]$ with arbitrary closed bounded intervals ; 

(b) $u$ (resp. $v$) of $\G$ with $h(u)$ (resp. $\ell(v)$) where $h$ and $\ell$ are continuous functions defined  on $U=[a,b]$ and $V=[c,d]$ respectively. In this case $U'=h(U)$ and $V'=\ell(V)$ are also bounded closed intervals. The Nash equilibrium point $(\bu,\bv)$ exists and is still of bang-bang type. The unique difference is that, it will jump between the  bound of set $U'$ (resp. $V'$) instead of $U$ (resp. $V$). \qedh
 \bigskip
 
\no \textbf{(ii)} As we indicated in Remark \ref{re:sigma}, the dynamics of the process $X^{t,x}$ of \eqref{eq: SDE} may contain a diffusion term $\sigma (t,x)$ (see equation \eqref{eq:sde sigma}) which is a function defined as:
  $$
(t,x)\in [0,T]\times \bR \mapsto \sigma(t,x) \in \bR,
  $$
with the following assumption:
\ms

\noindent \textbf{Assumption (A1)}: The function $\sigma(t,x)$ is uniformly Lipschitz w.r.t. $x$, it is invertible and bounded and its inverse is bounded.
\medskip

\no Under (A1), for any $\tx$, the following SDE 
  \begin{equation}\label{eq:sde sigma}\begin{array}{c}
  X_s^{t,x}=x+\int_t^s \sigma(r, X_r^{t,x})dB_r,\ \forall s\in [t,T]\text{ and } X_s^{t,x}=x \text{ for } s\in [0,t],\end{array}
  \end{equation}
has a unique solution (see e.g. Karatzas and Shreve, pp.289, \cite{karatzas}). Moreover $\sigma$ satisfies the uniform elliptic condition, \textit{i.e.} there exists a constant $\Upsilon>0$ such that for any $(t,x)\in [0,T]\times \bR$, $\Upsilon\leq \sigma(t,x)^2\leq \Upsilon^{-1}$. 
\ms

\no In this framework, the Hamiltonian functions associated with the NZSDG of payoffs given by (\ref{eq: payoffs}) are defined from $[0,T]\times \bR\times \bR \times U\times V$ into $\bR$ by:
  \begin{eqnarray*}
H_1(t,x,p,u,v)&:=p\sigma^{-1}(t,x)\Gamma(t,x,u,v)=p\sigma^{-1}(t,x)(f(t,x)+u+v);\\
H_2(t,x,q,u,v)&:=q\sigma^{-1}(t,x)\Gamma(t,x,u,v)=q\sigma^{-1}(t,x)(f(t,x)+u+v).
\end{eqnarray*}
Noticing that $\sigma^{-1}$ is bounded, it follows by the generalized Isaacs' condition (\ref{eq:Isaacs}) and the same approach in this article that, the Nash equilibrium point exists and is of bang-bang type. 

\no Actually we should point out that, all the results in this article will hold by the same techniques in this case (ii) with only some minor adaptions, except the convergence to $0$ of the second term on the right side of inequality \eqref{eq:eta in - eta im} which needs to be checked carefully. Indeed the objective is to show that 
   \begin{equation}\label{eq:sec term in gen}\ba{l}
|\E[\int_{t+\delta}^T (H_i^n(s, X_s^{t,x})-H_i^m(s,X_s^{t,x}))\cdot 1_{\{|X_s^{t,x}|\leq k\}}ds]\big|\rightarrow_{n,m\rightarrow \infty}0.\ea
   \end{equation}\medskip
for fixed $(t,x)$ and $k$. To begin with we give the following result related to domination of laws of the process $X^{t,x}$. 
   \begin{lemma}\label{lemma:domi}
{\bf ($L^q$-Domination)} Let $t\in [0,T]$ and $x,x_0\in \bR$. For 
$s\in [t,T]$, we denote by $\mu(t,x;s,dy)$ the law of $X_s^{t,x}$. Under Assumption (A1) on $\sigma$, for any $\bar q\in (1,\infty)$, the family of laws $\{\mu(t,x;s,dy),\ s\in (t,T]\}$ is $L^{\bar q}$-dominated by $\{\mu(0,x_0;s,dy),\ s\in (t,T]\}$, \textit{i.e.}, for any $\delta\in (0,T-t)$, there exists an application $\phi_{t,x,x_0}^{\delta}:[t+\delta, T]\times \bR\rightarrow \bR^+$ such that:
     
(a) $\mu(t,x;s,dy)ds=\phi_{t,x,x_0}^{\delta}(s,y)\mu(0,x_0;s,dy)ds$ for any $(s,y)\in [t+\delta,T]\times\bR ;$

(b) $\forall k\geq 1,\ \phi_{t,x,x_0}^{\delta}(s,y)\in L^{ \bar q}([t+\delta,T]\times [-k,k];\ \mu(0,x_0;s,dy)ds)$.
       \end{lemma}
   
\no \textit{Proof.} Readers are referred to \cite{H1997} (Section 28, pp.123) and \cite{H2014} (Lemma 4.3 and Corollary 4.4, pp.14-15) for the proof of this Lemma. However basically it uses the Aronson estimates \cite{aronson} for densities of the laws of the solution of SDE (\ref{eq:sde sigma}) under Assumption (A1). 
\bigskip

\no {\bf Proof of convergence \eqref{eq:sec term in gen}:} Thanks to Lemma \ref{lemma:domi}, there exists a function $\phi_{t,x,x_0}^{\delta}:[t+\delta,T]\times \bR\rightarrow \bR^+$ such that:
\begin{equation}\label{eq: ge phi}
\forall k\geq 1,\ \phi_{t,x,x_0}^{\delta}(s,y)\in L^{\frac{q}{q-1}}([t+\delta, T]\times [-k,k];\ \mu(0,x_0;s,dy)ds).
\end{equation}
Then 
   $$\ba{ll}
   |\E[&\!\!\!\!\!\int_{t+\delta}^T (H_i^n(s, X_s^{t,x})-H_i^m(s,X_s^{t,x}))\cdot 1_{\{|X_s^{t,x}|\leq k\}}ds]\big|\\\\
   {}&=|\int_{\bR}\int_{t+\delta}^T (H_i^n(s, y)-H_i^m(s,y))\cdot 1_{\{|y|\leq k\}}\mu(t,x;s,dy)ds|\\\\
  {} &=|\int_{\bR}\int_{t+\delta}^T (H_i^n(s, y)-H_i^m(s,y))\cdot 1_{\{|y|\leq k\}}\phi_{t,x,x_0}^{\delta}(s,y)\mu(0,x_0;s,dy)ds|.\ea
  $$
The constant $q$ in \eqref{eq: ge phi} is the one which makes that $H_i^n(s,y)\rightarrow_{n\rw \infty} H_i(s,y)$ weakly in $L^q([0,T]\times \bR;\ \mu(0,x_0;s,dy)ds)$ for $i=1,2$ and a fixed $q\in (1,2)$. Then combining this weak convergence result and \eqref{eq: ge phi} yields \eqref{eq:sec term in gen}.\qedh
\bigskip

\no \textbf{(iii)} In the same way one can deal with the multi-dimensional case for diffusion processes $X^{t,x}$ satisfying (\ref{eq:sde sigma}) when $\sigma (t,x)$ verifies Assumption (A1). \qedh

\end{document}